\documentclass[12pt]{article}

\usepackage{amssymb}
\usepackage{amsmath}

\newtheorem{thm}{Theorem}
\newtheorem{prop}{Proposition}
\newtheorem{defi}{Definition}
\newtheorem{pf}{Proof}

\newtheorem{rem}{Remark}

\newtheorem{lem}{Lemma}

\newcommand{\ai}{\operatorname{Ai}}
\newcommand{\ram}{\operatorname{A}}
\newcommand{\res}{\operatorname{Res}}

\begin{document}
\title{A connection formula between the Ramanujan function and the $q$-Airy function}
\setcounter{footnote}{2}
\author{Takeshi MORITA\thanks{Graduate School of Information Science and Technology, Osaka University,
1-1  Machikaneyama-machi, Toyonaka, 560-0043, Japan. 
E-mail:\ t-morita@cr.math.sci.osaka-u.ac.jp }}
\date{ }
\maketitle

\begin{abstract}
We show a connection formula between two different $q$-Airy functions. One is called the Ramanujan function which appears in Ramanujan's "Lost notebook". Another one is called the $q$-Airy function that obtained in the study of the second $q$-Painlev\'e equation. We use  the $q$-Borel transformation and the $q$-Laplace transformation following C. Zhang to obtain the connection formula.
\end{abstract}
%%%%%%%%%%%%%%%%%%%%%%%%%%%%%%section 1%%%%%%%%%%%%%%%%%%%%%%%%%%%%%%%%%%%%%%%
\section{Introduction}
In the study of the $q$-analysis, it is known that there exist several different $q$-special functions corresponding to a special function defined by a differential equation. For example, three types of $q$-Bessel functions are known. We denote them $J_\nu^{(1)}(x;q)$, $J_\nu^{(2)}(x;q)$ and $J_\nu^{(3)}(x;q)$ due to Ismail \cite{Is}. The first and the second one are called Jackson's first and second $q$-Bessel function and the third one is called the Hahn-Exton $q$-Bessel function. Similarly, two types of $q$-Airy functions are known. We denote them $\ram_q(x)$ and $\ai_q(x)$. The first one is called the Ramanujan function and the second one is called the $q$-Airy function. 

Three $q$-Bessel functions are given by 
\begin{align*}
J_\nu^{(1)}(x;q)&:=\frac{(q^{\nu +1};q)_\infty}{(q;q)_\infty}\left(\frac{x}{2}\right)^\nu\sum_{n\ge 0}\frac{1}{(q^{\nu +1};q)_n}\left(-\frac{x^2}{4}\right)^n,\\
J_\nu^{(2)}(x;q)&:=\frac{(q^{\nu +1};q)_\infty}{(q;q)_\infty}\left(\frac{x}{2}\right)^\nu\sum_{n\ge 0}\frac{q^{n^2}}{(q^{\nu +1};q)_n}\left(-\frac{q^\nu x^2}{4}\right)^n,\\
J_\nu^{(3)}(x;q)&:=\frac{(q^{\nu +1};q)_\infty}{(q;q)_\infty}x^\nu\sum_{n\ge 0}\frac{q^{\frac{n(n+1)}{2}}}{(q^{\nu +1};q)_n}\left(-x^2\right)^n.\\
\end{align*}

The Ramanujan function and the $q$-Airy function are defined by
\begin{align*}
\ram_q(x)&:=\sum_{n\ge 0}\frac{q^{n^2}}{(q;q)_n}(-x)^n,\\
\ai_q(x)&:=\sum_{n\ge 0}\frac{1}{(-q,q;q)_n}\left\{(-1)^nq^{\frac{n(n-1)}{2}}\right\}(-x)^n.
\end{align*}

Here, $(a;q)_n$ and $(a;q)_\infty$ are the $q$-Pochhammer symbol defined in section two. 
According to \cite{Andre} , the Ramanujan function $\ram_q(x)$ appears in the third identity on p.$57$ of Ramanujan's "Lost notebook" \cite{Ramanujan} as follows (with $x$ replaced by $q$):
\[\ram_q(-a)=\sum_{n\ge 0}\frac{a^nq^{n^2}}{(q;q)_n}=
\prod_{n\ge 1}\left(1+\frac{aq^{2n-1}}{1-q^ny_1-q^{2n}y_2-q^{3n}y_3-\cdots }\right)\]
where
\begin{align*}
y_1&=\frac{1}{(1-q)\psi^2(q)},\\
y_2&=0,\\
y_3&=\frac{q+q^3}{(1-q)(1-q^2)(1-q^3)\psi^2(q)}
-\frac{\sum_{n\ge 0}\frac{(2n+1)q^{2n+1}}{1-q^{2n+1}}}{(1-q)^3\psi^6(q)},\\
y_4&=y_1y_3,\\
\psi (q)&=\sum_{n\ge 0}q^{\frac{n(n+1)}{2}}=\frac{(q^2;q^2)_\infty}{(q;q^2)_\infty}.\\
\end{align*}

Ismail has pointed out that the Ramanujan function is a $q$-analogue of the Airy function \cite{Is}. 
We show more detail about these functions in section three.

It is known that there exist a relation for the $q$-Bessel functions and the $q$-Airy function. One is a relation between Jackson's first and second $q$-Bessel function \cite{Hahn} :
\begin{equation}
J_\nu^{(2)}(x;q)=\left(-\frac{x^2}{4};q\right)_\infty J_\nu^{(1)}(x;q).\label{ha}
\end{equation}
Another is a relation between the Hahn-Exton $q$-Bessel function and the $q$-Airy function \cite{Oyama} :
\begin{equation}
J_\nu^{(3)}(x;q)=\frac{(-q;q)_\infty}{(q;q)_\infty}x^\nu\ai_q(-qx^2),\label{oh}
\end{equation}
where $q^{\nu}=-1$.

Other relations are not known. The main result in this paper is a relation between two $q$-Airy functions. It is not known any relation between the Ramanujan function and the $q$-Airy function, but these two functions are related by a connection formula and not by algebraic relation like (\ref{ha}) or (\ref{oh}).

Connection problems of linear $q$-difference equations between the origin and the infinity are studied by G. D. Birkhoff \cite{Birkhoff}. 
Watson gave a connection formula for the basic hypergeometric equation ${}_2\varphi_1$ in $1910$ \cite{W}:

\begin{align*}
{}_2 \varphi_1\left(a,b;c;q,x \right)=&
\frac{(b,c/a;q)_\infty (a x,q/a x;q)_\infty }{(c, b/a ;q)_\infty (  x,q/x ;q)_\infty } 
{}_2 \varphi_1\left(a,aq/c;aq/b;q,cq/abx \right) \nonumber \\
&+ 
\frac{(a,c/b;q)_\infty (b x, q/ b x;q)_\infty }{(c, a/b;q)_\infty (  x,q/   x;q)_\infty } 
{}_2 \varphi_1\left(b, bq/c; bq/a; q, cq/abx \right). 
\end{align*}

Recently, C.~Zhang has given some connection formulae of $q$-difference equations of the confluent type \cite{Z0}, \cite{Z1} and \cite{Z2}. Zhang gives a connection formula of Jackson's $q$-Bessel function $J_\nu^{(1)}(x;q)$ \cite{Z2}. In \cite{Z2}, Zhang introduced the $q$-Borel transformation and the $q$-Laplace transformation which are useful to study connection problems. In section four, we apply Zhang's method to the $q$-Airy functions. 

The connection formula of the $q$-Airy function gives a relation between the Ramanujan function and the $q$-Airy function as follows:
\begin{thm}For any $x\in\mathbb{C}^*$,
\[\ram_{q^2}\left(-\frac{q^3}{x^2}\right)
=\frac{1}{(q,-1;q)_\infty}
\left\{\theta \left(\frac{x}{q}\right)\ai_q(-x)+\theta\left(-\frac{x}{q}\right)\ai_q(x)\right\}.\]
\end{thm}
Since our new relation shows an asymptotic behavior of the Ramanujan function near the infinity, it may be useful to study the Ramanujan function or similar type $q$-series. 
%%%%%%%%%%%%%%%%%%%%%%%%%%%section 2%%%%%%%%%%%%%%%%%%%%%%%%%%%%%%%%%%%%%%%%%%%
\section{Standard notations}
In this section, we fix our notations. We assume that $q\in\mathbb{C}^*$ satisfies $0<~|q|<1$. We define the $q$-Pochhammer symbol $(a;q)_n$. 
\begin{defi}For any $n\in\mathbb{Z}_{\ge 0}$,
\[(a;q)_n=
\begin{cases}
1,&\qquad n=0,\\
(1-a)(1-aq)\cdots (1-aq^{n-1}),&\qquad n\ge 1,
\end{cases}
\]
and 
\[(a;q)_\infty =\lim_{n\to\infty}(a;q)_n.\]
Moreover,
\[(a_1,a_2,\cdots ,a_m;q)_\infty =(a_1;q)_\infty (a_2;q)_\infty\cdots (a_m;q)_\infty .\]
\end{defi}

The $q$-difference operator $\sigma_q$ is given by $\sigma_qf(x)=f(qx)$. The basic hypergeometric series ${}_r\varphi_s$ is defined as follows.

\begin{defi}The basic hypergeometric series is given by
\[{_r\varphi_s}(a_1,\cdots ,a_r;b_1,\cdots ,b_s;q,x)
:=\sum_{n\ge 0}\frac{(a_1,\cdots ,a_r;q)_n}{(b_1,\cdots ,b_s;q)_n(q;q)_n}
\left[(-1)^nq^{\frac{n(n-1)}{2}}\right]^{1+s-r}x^n.\]
\end{defi}

We define the theta function of Jacobi. We denote by $\theta_q(x)$ or more shortly $\theta (x)$. The theta function of Jacobi is given by following series;
\begin{defi}For any $x\in\mathbb{C}^*$,
\[\theta_q(x)=\theta (x)=\sum_{n\in\mathbb{Z}}q^{\frac{n(n-1)}{2}}x^n.\]
\end{defi}
The theta function has some important properties. The following lemma is called Jacobi's triple product identity. 
\begin{lem}
For any $x\in\mathbb{C}^*$, we have
\[\theta (x)=\left(q,-x,-\frac{q}{x};q\right)_\infty .\]
\end{lem}
The theta function satisfies the following $q$-difference relation.
\begin{lem}\label{qt}
For any $k\in\mathbb{Z}, \theta (x)$ satisfies 
\[\theta (q^kx)=q^{-\frac{k(k-1)}{2}}x^{-k}\theta (x),\qquad \forall x\in\mathbb{C}^*.\]
\end{lem}
From lemma\ref{qt} , we remark that the function $\theta (-\lambda x)/\theta (\lambda x)$, $\forall\lambda\in\mathbb{C}^*$ satisfies a $q$-difference equation
\[u(qx)=-u(x)\]
which is also satisfied by the function $u(x)=e^{\pi i\left(\frac{\log x}{\log q}\right)}$. 
From the definition, the theta function has the following inversion formula.
\begin{lem}\label{inv}
For any $x\in\mathbb{C}^*$ , one gets
\[x\theta\left(\frac{1}{x}\right)=\theta (x).\]
\end{lem}

%%%%%%%%%%%%%%%%%%%%%%%%%%%%%section 3%%%%%%%%%%%%%%%%%%%%%%%%%%%%%%%%%%%%%%%

\section{Two types of the $q$-analogue of the Airy function}
There are two different $q$-analogue of the Airy function. One is called the Ramanujan function which appears in \cite{Ramanujan}. Ismail \cite{Is} pointed out that the Ramanujan function can be considered as a $q$-analogue of the Airy function. The other one is called the $q$-Airy function which is obtained by K.~ Kajiwara, T. Masuda, M. Noumi, Y. Ohta and Y. Yamada \cite{KMNOY} . In this section, we see the properties of these functions. We explain the reason why they are called $q$-analogue of the Airy function and we show $q$-difference equations which they satisfy.

%%%%%%%%%%%%%%%%%%%%%%%%%%%%%%%%%%section3.1%%%%%%%%%%%%%%%%%%%%%%%%%%%%%%%%%%%
\subsection{The Ramanujan function $\ram_q(x)$}
The Ramanujan function appears in Ramanujan's "Lost notebook" ~\cite{Ramanujan}.
Ismail has pointed out that the Ramanujan function can be considered as a $q$-analogue of the Airy function.
The Ramanujan function is defined by following convergent series;
\[\ram_q(x):=\sum_{n\ge 0}\frac{q^{n^2}}{(q;q)_n}(-x)^n
={_0\varphi_1}(-;0;q,-qx).\]

In the theory of ordinary differencial equations, the term Plancherel-Rotach asymptotics refers to asymptotics around the largest and smallest zeros. With $x=\sqrt{2n+1}-2^{\frac{1}{2}}3^{\frac{1}{3}}n^{\frac{1}{6}}t$ and for $t\in\mathbb{C}$, the Plancherel-Rotach asymptotic formula for Hermite polynomials $H_n(x)$ is
\begin{equation}
\lim_{n\to +\infty}\frac{e^{-\frac{x^2}{2}}}{3^{\frac{1}{3}}\pi^{-\frac{3}{4}}2^{\frac{n}{2}+\frac{1}{4}}\sqrt{n!}}H_n(x)=\ai (t). \label{pr}
\end{equation}
In \cite{Is}, Ismail shows the $q$-analogue of (\ref{pr});
\begin{prop}One can get 
\[\lim_{n\to\infty}\frac{q^{n^2}}{t^n}h_n(\sinh\xi_n|q)=\ram_q\left(\frac{1}{t^2}\right)\]
where $e^{\xi_n}=tq^{-\frac{n}{2}}$.
\end{prop}
Here, $h_n(\cdot |q)$ is the $q$-Hermite polynomial. 
In this sense, we can deal with the Ramanujan function $\ram_q(x)$ as a $q$-analogue
of the Airy function. The Ramanujan function satisfies the following $q$-diference equation;
\begin{equation}
\left(qx\sigma_q^2-\sigma_q+1\right)u(x)=0.
\label{ram}
\end{equation}

\begin{rem}We remark that another solution of the equation (\ref{ram}) is given by
\[u(x)=\theta (x){}_2\varphi_0(0,0;-;q,-x).\]
Here, 
\[{}_2\varphi_0(0,0;-;q,-x)=\sum_{n\ge 0}\frac{1}{(q;q)_n}\left\{(-1)^nq^{\frac{n(n-1)}{2}}\right\}^{-1}(-x)^n\]
is a divergent series.
\end{rem}

%%%%%%%%%%%%%%%%%%%%%%%%%%%%%%%%%%section3.2%%%%%%%%%%%%%%%%%%%%%%%%%%%%%%%%%%%
\subsection{The $q$-Airy function $\ai_q(x)$}
The $q$-Airy function is found by K. Kajiwara, T. Masuda, M. Noumi, Y. Ohta and Y. Yamada \cite{KMNOY}, in their study of the $q$-Painlev\'e equations. This function is the special solution of the second $q$-Painlev\'e equations and given by the following series
\[\ai_q(x):=\sum_{n\ge 0}\frac{1}{(-q,q;q)_n}\left\{(-1)^nq^\frac{n(n-1)}{2}\right\}(-x)^n={}_1\varphi_1(0;-q;q,-x).\]

T. Hamamoto, K. Kajiwara, N. S. Witte \cite{hama} proved following asymptotic expansions;

\begin{prop}
With $q=e^{-\frac{\delta^3}{2}}$, $x=-2ie^{-\frac{s}{2}\delta^2}$ as $\delta\to 0$,
\[{_1\varphi_1}(0;-q;q,-qx)=2\pi^{\frac{1}{2}}\delta^{-\frac{1}{2}}
e^{-\left(\frac{\pi i}{\delta^3}\right)\ln 2+\left(\frac{\pi i}{2\delta }\right)s+\frac{\pi i}{12}}\left[\ai\left(se^{\frac{\pi i}{3}}\right)+O(\delta^2)\right],\]
\[{_1\varphi_1}(0;-q;q,qx)=2\pi^{\frac{1}{2}}\delta^{-\frac{1}{2}}
e^{-\left(\frac{\pi i}{\delta^3}\right)\ln 2-\left(\frac{\pi i}{2\delta }\right)s-\frac{\pi i}{12}}\left[\ai\left(se^{-\frac{\pi i}{3}}\right)+O(\delta^2)\right]\]
for $s$ in any compact domain of $\mathbb{C}$.  
\end{prop}

Here, $\ai (\cdot )$ is the Airy function. From this proposition, we can regard the $q$-Airy function as a $q$-analogue of the Airy function.

We can easily check out that the $q$-Airy function satisfies the second order linear $q$-difference equation
\begin{equation}
\left(\sigma_q^2+x\sigma_q-1\right)u(x)=0.
\label{qai}
\end{equation}
Another solution of the equation (\ref{qai}) is given by
\[u(x)=
e^{\pi i\left(\frac{\log x}{\log q}\right)}{_1\varphi_1}(0;-q;q,x)
=e^{\pi i\left(\frac{\log x}{\log q}\right)}\ai_q(-x).\]

%%%%%%%%%%%%%%%%%%%%%%%%%%%%%%%%%%section3.3%%%%%%%%%%%%%%%%%%%%%%%%%%%%%%%%%%%%

\subsection{Shearing transformations}
We define a shearing transformation of a second order linear $q$-difference equation.
\begin{defi}For a $q$-difference equation
\begin{equation}\label{sh}
a(x)u(q^2x)+b(x)u(qx)+c(x)u(x)=0,
\end{equation}
we define the shearing transformation as follows
\[t^2:=x,\quad v(t):=u(t^2),\quad p:=\sqrt{q}.\]
\end{defi}
The shearing transform of the equation (\ref{sh}) is given by
\[a(t^2)v(p^2t)+b(t^2)v(pt)+c(t^2)v(t)=0.\]
By the shearing transformation, the equation 
\[\left(K\cdot x\sigma_q^2-\sigma_q+1\right)u(x)=0\]
is transformed to
\begin{equation}
\left(K\cdot t^2\sigma_p^2-\sigma_p+1\right)v(t)=0,\label{kram}
\end{equation}
where $K$ is a fixed constant in $\mathbb{C}^*$.

%%%%%%%%%%%%%%%%%%%%%%%%%%%%%%%%%%section3.4%%%%%%%%%%%%%%%%%%%%%%%%%%%%%%%%%%%%
\subsection{The $q$-Airy equation around the infinity}
We consider the behavior of the equation (\ref{qai}) around the infinity. We set $x=1/t$ and $z(t)=u(1/t)$. Then $z(t)$ satisfies 
\[\left(-\sigma_q^2+\frac{1}{q^2t}\sigma_q+1\right)z(t)=0.\]
We set $\mathcal{E}(t)=1/\theta (-q^2t)$ and $f(t)=\sum_{n\ge 0}a_nt^n,\quad a_0=1$. We assume that $z(t)$ can be described as 
\[z(t)=\mathcal{E}(t)f(t)=\frac{1}{\theta (-q^2t)}\left(\sum_{n\ge 0}a_nt^n\right).\]

The function $\mathcal{E}(t)$ has the following property;
\begin{lem} For any $t\in\mathbb{C}^*$,
\[\sigma_q\mathcal{E}(t)=-q^2t\mathcal{E}(t),
\qquad \sigma_q^2\mathcal{E}(t)=q^5t^2\mathcal{E}(t).\]
\end{lem}
From this lemma, $f(t)$ satisfies the following equation
\begin{equation}\label{qaeq}
\left(-q^5t^2\sigma_q^2-\sigma_q+1\right)f(t)=0.
\end{equation}
Since (\ref{qaeq}) is the same as (\ref{kram}) for $K=-q^5$, we obtain
\[f(t)={}_0\varphi_1(-;0;q^2,q^5t^2)=\ram_{q^2}(-q^3t^2).\]

%%%%%%%%%%%%%%%%%%%%%%%%%%%%section 4%%%%%%%%%%%%%%%%%%%%%%%%%%%%%%%%%%%%%%%%%

\section{The $q$-Borel transformation, the $q$-Laplace transformation and the connection formula}
In this section, we show a connection formula for $f(t)$. In order to obtain a connection formula, we need the $q$-Borel transformation and the $q$-Laplace transformation following Zhang \cite{Z1}.
\subsection{\label{qbql}The $q$-Borel transformation and the $q$-Laplace transformation}
\begin{defi}For $f(t)=\sum_{n\ge 0}a_nt^n$, the $q$-Borel transformation is defined by
\[g(\tau )=\left(\mathcal{B}_qf\right)(\tau ):=\sum_{n\ge 0}a_nq^{-\frac{n(n-1)}{2}}\tau^n,\]
and the $q$-Laplace transformation is given by
\[\left(\mathcal{L}_qg\right)(t):=\frac{1}{2\pi i}\int_{|\tau |=r}g(\tau )\theta\left(\frac{t}{\tau}\right)\frac{d\tau}{\tau },\qquad 0<r<\frac{1}{|q^2|}.\] 
\end{defi}
 The $q$-Borel transformation can be considered as a formal inverse of the $q$-Laplace transformation.
\begin{lem}For any entire function $f$,
\[\mathcal{L}_q\circ\mathcal{B}_qf=f.\]
\end{lem}
\begin{pf}We can prove this lemma calculating residues of the $q$-Laplace transformation around the origin. \rule[-2pt]{5pt}{10pt}
\end{pf}

The $q$-Borel transformation has following operational relation;
\begin{lem}\label{orqb}
For any $l,m\in\mathbb{Z}_{\ge 0}$,
\[\mathcal{B}_q(t^m\sigma_q^l)=q^{-\frac{m(m-1)}{2}}\tau^m\sigma_q^{l-m}\mathcal{B}_q.\]
\end{lem}
\subsection{The connection formula of the $q$-Airy function}
Applying the $q$-Borel transformation in \ref{qbql} to the equation (\ref{kram}) and using lemma \ref{orqb}, we obtain the first order $q$-difference equation
\[g(q\tau )=(1+q^2\tau )(1-q^2\tau )g(\tau ).\]
Since $g(0)=1$, $g(\tau )$ is given by an infinite product
\[g(\tau )=\frac{1}{(-q^2\tau ;q)_\infty(q^2\tau ;q)_\infty}\]
which has single poles at
\[\left\{\tau ;\tau =\pm q^{-2-k},\quad \forall k\in\mathbb{Z}_{\ge 0}\right\}.\]
By Cauchy's residue theorem, the $q$-Laplace transform of $g(\tau )$ is 
\begin{align*}
f(t)=&\frac{1}{2\pi i}\int_{|\tau |=r}g(\tau )\theta\left(\frac{t}{\tau}\right)\frac{d\tau}{\tau }\\
=&-\sum_{k\ge 0}\res\left\{g(\tau )\theta\left(\frac{t}{\tau}\right)\frac{1}{\tau };\tau =-q^{-2-k}\right\}\\
&-\sum_{k\ge 0}\res\left\{g(\tau )\theta\left(\frac{t}{\tau}\right)\frac{1}{\tau };\tau =q^{-2-k}\right\}
\end{align*}
where $0<r<r_0:=1/|q^2|$. We can culculate the residue from lemma \ref{lems} and lemma \ref{qt} .
\begin{lem}\label{lems}For any $k\in\mathbb{N}$, $\lambda\in\mathbb{C}^*$, one can get;
\begin{enumerate}
%\item $\theta_q(q^nt)=q^{-\frac{n(n-1)}{2}}t^{-n}\theta_q(t)$,
\item $\res\left\{\dfrac{1}{\left(\tau /\lambda ;q\right)_\infty}\dfrac{1}{\tau}:\tau =\lambda q^{-k}\right\}
=\dfrac{(-1)^{k+1}q^{\frac{k(k+1)}{2}}}{(q;q)_k (q;q)_\infty}$,
\item $\dfrac{1}{(\lambda q^{-k};q)_\infty}=\dfrac{(-\lambda )^{-k}q^{\frac{k(k+1)}{2}}}{(\lambda ;q)_\infty \left(q/\lambda ;q\right)_k},\quad \lambda \not \in q^{\mathbb{Z}}$.
\end{enumerate}
\end{lem}
Summing up all of residues, we obtain 
\begin{align*}
f(t)&=\frac{\theta (q^2t)}{(q,-1;q)_\infty} {_1\varphi_1}\left(0,-q;q,\frac{1}{t}\right)
+\frac{\theta (-q^2 t)}{(q,-1;q)_\infty} {_1\varphi_1}\left(0,-q;q,-\frac{1}{t}\right).
\end{align*}
Combining with lemma\ref{inv}, we get a connection formula for $z(t)=\mathcal{E}(t)f(t)$. Finally, we acquire the following connection formula between the Ramanujan function and the $q$-Airy function.
\begin{thm}
For any $x\in\mathbb{C}^*$,
\[\ram_{q^2}\left(-\frac{q^3}{x^2}\right)
=\frac{1}{(q,-1;q)_\infty}
\left\{\theta \left(\frac{x}{q}\right)\ai_q(-x)+\theta\left(-\frac{x}{q}\right)\ai_q(x)\right\}.\]
\end{thm}
Here, both $\ram_q(x)$ and $\ai_q(x)$ are defined by convergent series on whole of the complex plain. The connection formula above is valid for any $x\in\mathbb{C}^*$.
\section*{Acknowledgement}
The author express his thanks to Professor Yousuke Ohyama for his careful conduct and kindly encouragement. The author also expresses his thanks to Mr. Nobutaka~Nakazono from Kyushu University for his valuable comments. This work is partially supported by the Mitsubishi foundation. %% and the JSPS Grant-in-Aid for Scientific Research.


\begin{thebibliography}{99}

\bibitem{Andre}
G.~E.~Andrews, {\it Ramanujan's "lost" notebook. VIII: the entire Rogers-Ramanujan function,} Adv. Math.{\bf 191} (2005), no. 2, 393-407.




\bibitem{Birkhoff}
G.~D.~Birkhoff, The generalized Riemann problem for linear differential
equations and the allied problems for linear difference and $q$-difference
equations, Proc. Am. Acad. Arts and Sciences, $49$ ($1914$), $521-568$.


\bibitem{GR}
G.~Gasper and M.~Rahman (1990). Basic Hypergeometric Series . 
Encycl. Math. Appl. Cambridge Univ. Press , Cambridge. 


\bibitem {Hahn}
W.~Hahn, 
Beitr\"age zur Theorie der Heineschen Reihen. Die $24$ Integrale der Hypergeometrischen $q$-Differenzengleichung. 
Das $q$-Analogon der Laplace-Transformation. 
Math. Nachr.  {\bf 2},  (1949). 340--379. 


\bibitem{hama}
T.~Hamamoto, K.~Kajiwara, N.~S.~Witte, 
Hypergeometric solutions to the $q$-Painlev\'e equation of type 
$(A_1+A_1')^{(1)}$, Int. Math. Res. Not. Vol. 2006, Article ID 84619, Pages $1-26$.

\bibitem{Is}
M.~E.~H.~Ismail,
 Asymptotics of $q$-Orthogonal Polynomials
and a $q$-Airy Function, Int. Math. Res. Not. (2005), No. { 18} 1063--1088.

\bibitem{KMNOY}
K.~Kajiwara, T.~Masuda, M.~Noumi, Y.~Ohta, and Y.~Yamada, 
Hypergeometric solutions to the $q$-Painlev\'e equations, Int. Math. Res. Not.  (2004), no. 47, 2497--2521.

\bibitem{Oyama}
Y.~Ohyama 
A unified approach to $q$-special functions of the Laplace type,
{\ttfamily arXiv:1103.5232}.

\bibitem{Ramanujan}
S.~Ramanujan, The Lost Notebook and Other Unpublished Papers (Intro by G. ~E. ~Andrews), Narosa, New Delhi, 1988.

\bibitem{Z3}
R.~F.~Swarttouw and H.~G.~ Meijer
A $q$-analogue of the Wronskian and a second solution of
the Hahn-Exton $q$-Bessel difference equation
 Proc. Am. Math. Soc. 129 (1994),   855--864.

\bibitem{W}
G.~N.~Watson,  The continuation of functions defined by generalized
hypergeometric series, \textit{Trans. Camb. Phil. Soc.} {\bf 21} (1910),  281--299.

\bibitem{Z0}
C.~Zhang, Remarks on some basic hypergeometric series, 
in ``Theory and Applications of Special Functions", 
Springer (2005), 479--491. 

\bibitem{Z1}
C.~Zhang, Sur les fonctions $q$-Bessel de Jackson, 
  J. Approx. Theory, {\bf 122}  (2003),  208--223.  

\bibitem{Z2}
C.~Zhang, Une sommation discr\`e pour des \'equations aux $q$-diff\'erences lin\'eaires et \`a  coefficients, analytiques: th\'eorie g\'en\'erale et exemples, in ``Differential Equations and Stokes Phenomenon'', World Scientific (2002), 309--329.

\end{thebibliography}
\end{document}